\documentclass{IEEEcsmag}

\usepackage[colorlinks,urlcolor=blue,linkcolor=blue,citecolor=blue]{hyperref}
\expandafter\def\expandafter\UrlBreaks\expandafter{\UrlBreaks\do\/\do\*\do\-\do\~\do\'\do\"\do\-}
\usepackage{upmath,color}

\jvol{XX}
\jnum{XX}
\paper{8}
\jmonth{Month}
\jname{Publication Name}
\jtitle{Publication Title}
\pubyear{2023}

\usepackage{glossaries}
\newacronym{amr}{AMR}{adaptive mesh refinement}
\newacronym{cfd}{CFD}{computational fluid dynamics}
\newacronym{gpu}{GPU}{graphics processing unit}
\newacronym{cpu}{CPU}{central processing unit}
\newacronym{hpc}{HPC}{high performance computing}
\newacronym{fpga}{FPGA}{field-programmable gate array}
\newacronym{asic}{ASIC}{application-specific integrated circuit}
\newacronym{pde}{PDE}{partial differential equation}
\newacronym{mol}{MOL}{method of lines}
\newacronym{nrel}{NREL}{National Renewable Energy Laboratory}
\newacronym{doe}{DOE}{Department of Energy}
\newacronym[longplural={Gaussian processes}]{gp}{GP}{Gaussian Process}
\newacronym{ei}{EI}{expected improvement}
\newacronym{ac}{AC}{Adaptive Computing}
\newacronym{smt}{SMT}{Surrogate Modeling Toolbox}
\newacronym{uq}{UQ}{uncertainty quantification}
\newacronym{mme}{MME}{metal modulated epitaxy}
\newacronym{md}{MD}{molecular dynamics}
\newacronym{amar}{AMAR}{Adaptive Mesh and Algorithm Refinement}
\newacronym{kmc}{KMC}{kinetic Monte-Carlo}
\newacronym{mpi}{MPI}{Message Passing Interface}
\newacronym{qoi}{QOI}{quantity of interest}
\newacronym{nn}{NN}{neural network}
\newacronym{hero}{HERO}{Hybrid Environment Resources and Operations}
\newacronym{rom}{ROM}{reduced-order model}
\newacronym{mf}{MF}{multi-fidelity}
\newacronym{ai}{AI}{artificial intelligence}
\newacronym{lpv}{LPV}{linear parameter varying}
\newacronym{bo}{BO}{Bayesian optimization}
\newacronym{lhs}{LHS}{Latin hypercube sampling}
\newacronym{eh}{EH}{enzymatic hydrolysis}
\newacronym{gan}{GaN}{gallium nitride}
\newacronym{mpc}{MPC}{model predictive control}

\makeglossaries
\glsdisablehyper
\glsenableentrycount

\usepackage{algorithm}
\usepackage{algorithmic}
\usepackage{amsmath}

\begin{document}

\sptitle{Theme Article: NREL Special Issue}

\title{Adaptive Computing for Scale-up Problems}

\author{Kevin Patrick Griffin*\thanks{* These authors contributed equally}, Hilary Egan*, Marc T. Henry de Frahan, Juliane Mueller, Deepthi Vaidhynatha, Dylan Wald, Rohit Chintala, Olga A. Doronina, Hariswaran Sitaraman, Ethan Young, Ryan King, Jibonananda Sanyal, Marc Day}
\affil{National Renewable Energy Laboratory, Golden, CO, 80401, USA}

\author{Ross E. Larsen}
\affil{National Renewable Energy Laboratory, Golden, CO, 80401, USA and Renewable and Sustainable Energy Institute, University of Colorado Boulder, Boulder, Colorado 80309}

\begin{abstract}
Adaptive Computing is an application-agnostic outer loop framework to strategically deploy simulations and experiments to guide decision making for scale-up analysis. Resources are allocated over successive batches, which makes the allocation adaptive to some objective such as optimization or model training. The framework enables the characterization and management of uncertainties associated with predictive models of complex systems when scale-up questions lead to significant model extrapolation. A \textcolor{black}{key advancement of this framework is its integration of multi-fidelity surrogate modeling, uncertainty management, and automated orchestration of various computing and experimentation resources into a single integrated software package. This enables efficient multi-fidelity modeling across multiple computing resources by incorporating real-world constraints such as relative queue times and throughput on individual machines into the multi-fidelity sampling decision.} We discuss applications of this framework to problems in the renewable energy space, including biofuels production, material synthesis, perovskite crystal growth, and building electrical loads. 
\end{abstract}


\maketitle
\chapteri{S}cale-up problems can be defined as the challenge of (1) translating a laboratory-scale experiment or small-scale simulation to a large-scale (e.g., industrial) application and (2) characterizing and managing the uncertainty introduced by scaling up the devices or processes involved thus leading to extrapolation beyond known data regimes. It is often inaccurate to directly extrapolate the results of laboratory-scale experiments to large-scale systems because new length and time scales can alter the balance of physical processes. For example, longer mixing times in a full-scale bioreactor might render it uneconomical at scale. 

Computer simulations are essential to address issues of scale-up since they allow for direct evaluation of physical processes known to exist at the larger scales. However, computational costs typically increase quickly with growing length and time scales, and, in many applications, only a limited number of full-scale calculations are affordable, even with ever-growing computational resources. Engineering design optimization typically requires a large number of simulated model realizations, and this ultimately limits the fidelity, and associated cost, of the underlying models that can be employed.
This has motivated the expanding use of multi-fidelity modeling strategies, which leverage a large number of lower-fidelity (lower-cost, lower-accuracy) calculations to learn trends in the data, with a limited number of expensive higher-fidelity evaluations (from simulations or experiments) that inform corrections to the lower-fidelity predictions. 

Multi-fidelity modeling orchestrates the use of multiple models and/or experiments to efficiently characterize key underlying processes. Though these specialized modeling strategies can be used to capture distinct regimes of interest, they have complex convergence properties, uncertainties that are hard to quantify, or stochasticities that make their use considerably more challenging than single-model systems. Furthermore, while many ad-hoc methods and specialized analyses have been developed for comparison of disparate, often multi-modal, data sources, these works typically assume that the sample space is identical across fidelity levels. Meanwhile, for problems that relate to scale-up, the sample space typically varies across the fidelity levels. For example, laboratory experiments cannot be run at the same operating conditions as a full-scale reactor, where simulations may be especially useful.

Beyond differences in the modalities and representations of the underlying models, the drastic increase in the heterogeneity of computing resources has driven the development of a similar diversity in specialized workflows. These trends include \gls{cpu}-based \gls{hpc} resources being supplanted by heterogenous architectures including \glspl{gpu}, \glspl{fpga} and \glspl{asic}, \gls{hpc} centers being supplemented by flexible on-demand compute such as cloud-based services, and the proliferation of sensors and connected devices leading to increased opportunities for edge computing, i.e., discrete computing devices co-located with sensors to pre-process data and reduce data-load. \textcolor{black}{This diversity in computing architectures and workflows often lends itself to application-specific software for surrogate modeling and uncertainty management, limiting software reusability.} Similarly, there are many real-time or autonomous user facilities, laboratories, and other data generating resources across the \gls{doe}, which bring their own challenges with regard to data management, resource management, cross comparison of results, and making decisions within a limited time frame (especially when steering or guiding experiment with computation). 

In this paper we detail a computational outer-loop framework called \emph{\gls{ac}} to support decision making for scale-up analysis, where data acquisition is constrained to a specified budget. \textcolor{black}{\gls{ac} integrates} multi-fidelity surrogate modeling, \textcolor{black}{adaptive sampling, and automated resource management to coordinate computational and experimental resources for applications in optimization, multiscale modeling, and control. The result is a flexible and modular tool that accommodates real-world \gls{hpc} constraints, including resource schedulers, queue wait times, and heterogeneous computing environments. This integration enables more effective decision-making for large-scale simulations, ensuring that computational resources are allocated adaptively based on both model fidelity and system and resources availability.}
\section{AC Framework}
\begin{figure*}[hbt!]
\centering
\includegraphics[width=1.0\textwidth]{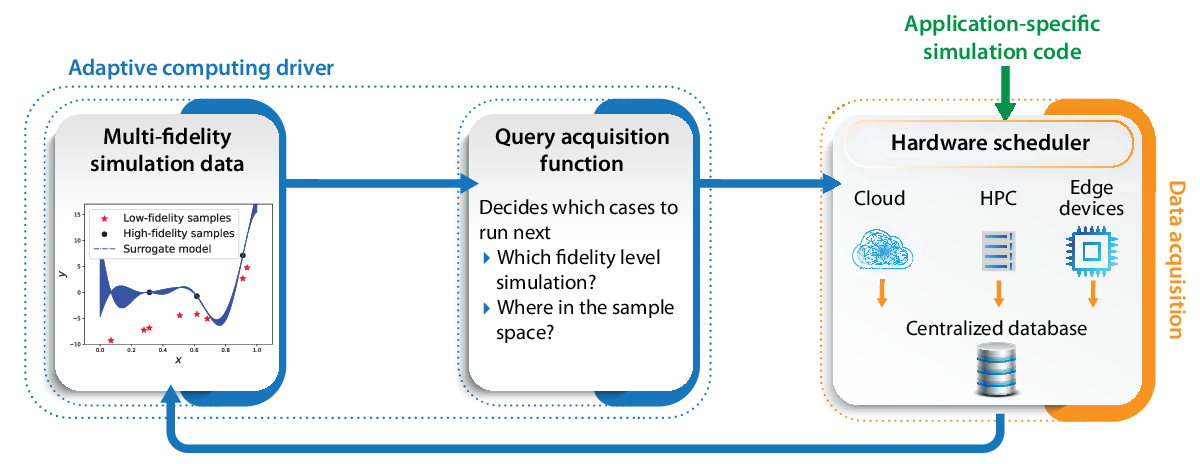}
\caption{The \gls{ac} software drives the scheduling of application-specific simulations. A surrogate model informs an acquisition function, which selects simulation cases. The hardware scheduler manages their execution on computational resources.}
\label{fig:ac_framework}
\end{figure*}

The fundamental tenet behind \gls{ac} is defining a scale-up problem through the lens of a goal-oriented outer-loop application. This outer-loop application could be an optimization, uncertainty quantification or sensitivity analysis task, all of which typically require large numbers of data samples to be collected and are thus challenging in settings with limited budgets. To address this challenge, the outer-loop \gls{ac} driver uses surrogate models and adaptive sampling strategies to strategically choose which data samples to acquire next using the most suitable type of simulation or experiment. The requested data samples are asynchronously evaluated and fed back into the \gls{ac} driver for informing the selection of the next batch of data samples to be acquired. This workflow is illustrated in Figure~\ref{fig:ac_framework}. 
Each step of the workflow can be accomplished using various methods. For example, the surrogate model can be a multi-fidelity deep \gls{nn} or \gls{gp} model. The total budget of the different types of acquirable samples can be fixed or dynamic. The method for choosing new data samples and the corresponding evaluation mechanism can rely on predictions from the surrogate model, uncertainty estimates, levels of trust, distance metrics, or random sampling. A key advantage of the \gls{ac} framework is the modular design of this infrastructure: each step can be tailored to the specific problem and they are tightly coupled to enable seamless computation across the feedback loop.

\subsection{Multi-Fidelity Surrogates} \label{sec:bo}

Multi-fidelity surrogate modeling is a high-level strategy for integrating data from separate sources (simulations, models, and/or experiments) into a single surrogate for comprehensive evaluation. Key considerations in the design of the global \gls{mf} surrogate model include: 1) the choice of underlying form of the surrogate 2) the design of the bridging function that relates model fidelities, and 3) the methods for communicating and propagating uncertainty. 

For a given problem the optimal surrogate model choice depends on the absolute cost of data acquisition and expected smoothness of the response. While \glspl{gp} are robust with limited training data, their expressiveness is limited for complex or non-smooth mappings; in contrast \glspl{nn} are much more data-hungry but can be universal function approximators. \textcolor{black}{There are also projection-based (intrusive/non-intrusive) model reduction methods such as the proper orthogonal decomposition. See \cite{peherstorfer2018survey} for several recent developments. Model reduction methods can be efficient for representing high-dimensional data, but these methods are typically best suited for data that is structured in the input space. Also, these methods typically do not inherently provide uncertainty estimates. In this work, we focus on more general input-output black-box machine-learning techniques that are sufficiently robust to support lower dimensional quantities of interest with sparse data and that provide a mechanism for uncertainty quantification.} 

In a multifidelity context, bridge functions can be used to connect the individual models across the fidelity hierarchy. 
A form for such a bridge function commonly used with GPs is an autoregressive model of order 1 \cite{Kennedy2001}:
\begin{equation} \label{eq:mf}
    y_{MF}[x] = \rho[x] y_{LF}[x] + \delta[x], 
\end{equation}
where $x$ is the sample point, $y_{LF}$ are low fidelity function values, $\rho$ is the scaling function, and $\delta$ is the discrepancy function. For example, $\rho$ and $\delta$ can be represented by polynomials with constants determined by least-squares regression of the multi-fidelity model $y_{MF}$ to high fidelity function values. \glspl{nn} can also be used for the bridging functions, and have similar advantages/disadvantages in smoothness and data requirements compared to a pure autoregressive approach here, but have additionally been useful for embedding higher level physical constraints \cite{meng2019}.

\subsection{Uncertainty and Trust Modeling}
Uncertainty management begins with accurate quantification and propagation of epistemic and aleatoric uncertainties \cite{ghanem2017handbook}. Epistemic uncertainty can arise both from the individual models in a fidelity hierarchy (due to undersampling) and from the bridging functions that connect them (e.g., due to model bias). When using \glspl{gp} as the underlying surrogate model, epistemic uncertainty is automatically quantified through samples in the latent Gaussian space; similarly other surrogate model classes can be used to manage epistemic uncertainty directly through, e.g., ensemble-based approaches. Aleatoric uncertainty can arise from inherent stochasticity in a simulation or experimental irrepeatability due to imprecise measurements or uncontrolled disturbances. 

For many observable quantities, uncertainty is normally distributed, and can be homoskedastic (constant), but it may also be heteroskedastic (varying) across the input parameter space. In either case, with sufficient data, it is possible to estimate uncertainty using hierarchical optimization of the parameters of the multi-fidelity surrogate. However, in problems that involve extrapolation, these uncertainty estimates usually explode or become artificially constrained by specified uncertainty priors, which are typically ad hoc and homoskedastic. Thus, these uncertainty estimates quickly become uninformative for data acquisition decisions related to scale-up. 

It is often useful to leverage domain expertise to provide additional information about the trust one has about the validity of different models for specific regimes. For example, the high-fidelity model in the hierarchy for a fluid system could be based on a continuum representation operating well for low Knudsen numbers $Kn$, the non-dimensional number which tends toward zero in the continuum limit. Errors for the flow rate in internal flows are expected to scale with $c_1+c_2Kn+c_3Kn^2$ \cite{hadjiconstantinou2006limits}, where $c_{\cdot}$ are model constants. The user can provide this known deficiency of the Navier-Stokes for weakly rarefied flows by specifying this error scaling in the uncertainty prior for all continuum models in the multi-fidelity modeling hierarchy. 

This model-specific local trust estimate is called a domain-expert-informed prior \cite{ghanem2017handbook} because it leverages physical insight (domain expertise) about the trustworthiness of a specific type of modeling assumption and it is a function of where the model is queried in the sample space, e.g., at which value of $Kn$. Similar trust estimates can be made for some experimental measurement techniques. For numerical methods, trust can depend on spatial and temporal resolutions (i.e., model hyperparameters); it can be smoothly varying, or discontinuous due, for example, to hidden constraints that prevent data acquisition in certain regions of parameter space.
This trust model can then be exploited, not only to guide data acquisition, but to provide essential information for decision making given the outputs obtained from \gls{ac}.

\subsection{Budget-Constrained Data Acquisition Strategies}
Given the multi-fidelity surrogate with quantified uncertainty, we define a sampling strategy to identify new candidate function evaluations in the sample space to achieve the outer-loop goal. The first ingredient is selecting a criterion for the benefit $r_{m\ell}$ of computing a simulation, which is \textcolor{black}{the} acquisition function \textcolor{black}{value divided by the cost of the simulation}. For example, for optimization problems, the \gls{ei} acquisition function is commonly used with \glspl{gp} \cite{Jones1998}. This function balances local and global search by combining predicted function values with their corresponding uncertainty and trust estimates provided by the GP. Other acquisition functions exist such as probability of improvement, optimizing a confidence bound, using weighted sums of predicted function values and distance to already evaluated points, etc. Each has advantages and disadvantages~\cite{shoemaker-paper}. 

In many multi-fidelity applications, we not only choose where to sample but also which fidelity level to evaluate while also optimizing over limited available resources. Moreover, sampling decisions may also be constrained by feasibility of an experiment or simulation at a given point in parameter space. AC's sampling strategy thus exploits predictions from the multi-fidelity surrogate, uncertainty and trust estimates, and fidelity-specific feasibility constraints to decide the most promising points in the parameter space, and an optimization over the available resources allows optimal distribution of the evaluations. The goal is to arrive at the best possible solutions within the limited budget and to provide uncertainty estimates together with the solution to enable decision making. \textcolor{black}{The resource budget and wall-clock limit per iteration, $B_i$ and $T_i$, respectively, can be set as constants or estimated dynamically using look-ahead strategies to choose the batch size using estimates of future needs and the remaining budget \cite{lam2016bayesian}}

An essential piece of the sampling strategy is batching, which permits parallel function evaluations \cite{shahriari2015taking}. Here, the sampling strategy determines the next candidate for evaluation. Rather than immediately conducting this evaluation, it is added to a queue of points to be evaluated. The expected value from the surrogate model is treated as a placeholder value for the corresponding result, and the surrogate is retrained using this value. This algorithm is applied successively to develop a batch of $M_\ell$ candidate input locations for each fidelity level $\ell$. There is a trade-off between the benefits of frequently updating the inputs to the sampling strategy and the wall-clock efficiency of parallel sample evaluation that must be factored into the ideal batch size for a given application.

Next, the optimization problem in Algorithm \ref{alg:mf_opt} is solved to decide which of the $M_\ell$ candidate input locations, i.e., which simulations or computer-controlled experiments, will be executed. We optimize for binary decision variables $y_{m\ell}$ with $y_{m\ell}=1$ if the $m^{\rm{th}}$ sample point is evaluated on fidelity level $\ell$ so as to maximize the cumulative benefit $\sum_{\ell=1}^L \sum_{m=1}^{M_\ell} r_{m\ell} y_{m\ell}$.
This maximization is subject to the constraints of the available budgets for the wall-clock time or resource utilization, which might be distinct for each fidelity level or shared across levels. Details for specific applications are discussed in the following sections.

\begin{algorithm}
\caption{Multi-Fidelity Data Acquisition Strategy}\label{alg:mf_opt}
    \begin{algorithmic}[1]
        \STATE \textbf{Input:}
        \STATE $L$ - Number of fidelity levels
        \STATE $I$ - Number of batches or iterations
        \STATE $T^r$ - Remaining wall-clock budget
        \STATE $B^r$ - Remaining resource budget
        
        \STATE \textbf{Output:}
        \STATE $y_{m\ell}$ - Values of binary decision variables for each candidate simulation $m$ on level $\ell$
        
        \FOR{$i = 1$ to $I$}
            \STATE Assign $T_i$ - Wall-clock limit for iteration $i$
            \STATE Assign $B_i$ - Resource budget for iteration $i$
            \STATE Choose $M_\ell$ - Number of input locations to consider
            \STATE Compute $t_{m\ell}$ - Computational cost of conducting a simulation
            \STATE Compute $r_{m\ell}$ - Benefit of evaluating a simulation
            \STATE Solve discrete optimization problem for $y_{m\ell}$ that maximizes the cumulative benefit subject to constraints $T_i$ and $B_i$
            \STATE $T^r \leftarrow T^r-T_i$
            \STATE $B^r \leftarrow B^r-B_i$
            \IF{$T^r<0$ \OR $B^r<0$}
                \STATE \textbf{break}	 			
            \ENDIF
        \ENDFOR
    \end{algorithmic}
\end{algorithm}

\subsection{Resource Management}
The solution to the previous optimization problem is a set of cases (i.e., tasks) to run on computational or experimental hardware. (See the ``data acquisition" step of Figure \ref{fig:ac_framework}.) The orchestration of these tasks and their data is handled using the \gls{hero} framework, which is a scalable, low-latency, and distributed task queue, capable of interfacing with the large number of heterogenous resources that may be necessary for a given problem (e.g., \gls{hpc} resources with or without \glspl{gpu}, a distributed cloud environment, edge devices such as network-attached microcontrollers, or experiment controllers). This scheme allows asynchronous execution across resources with vastly different time-scales and latencies, while simplifying the execution of on-demand model re-training based on non-blocking data availability.

The \gls{ac} driver program (see Figure \ref{fig:ac_framework}) adds tasks to \gls{hero} queues based on their resource requirements. Compute resources run workers, i.e., programs that check the queues for jobs that they can complete. The number of workers assigned to a queue determines the maximum throughput. The worker claims the task, executes the task, and records the simulation output in a centralized database. If applicable, workers can service multiple queues and can prioritize them. This system is resilient and robust for simultaneous queries and device failures which are more likely to occur with vast numbers of disparate resources. 


\section{Motivating Science Applications}
We demonstrate how the \gls{ac} framework is applied to a diverse array of problem classes, with a specific focus on recasting the illustrative renewable energy examples into the \gls{ac} workflow.

\subsection{Engineering Design Optimization: Biofuels Virtual Engineering}

An engineering design problem can often be framed as an optimization over allowable ranges of design parameters such as material compositions and geometrical specifications.
In the context of scale-up, design optimization has additional uncertainties since design decision must be made from a limited number of simulations of large-scale systems and many small-scale simulations and experiments.
For illustration, we consider design optimization for an economically viable design of a chemical reactor for processing second-generation biofuels, which can be generated from waste (\textcolor{black}{lignocellulosic}) biomass such as the non-edible parts of crops.

\subsubsection{AC formulation:}
The \gls{ac} outer loop for this application is the design optimization problem. \gls{ac} selects the input values (design specifications) to be simulated next based on surrogate-based optimization of some design objective. 
Specifically, \gls{eh} simulations are conducted using an \gls{nrel} software package called \textcolor{black}{\textit{Vebio}, a Virtual engineering framework with several zero and multidimensional models for biomass conversion processes \cite{Young2023}}. \textcolor{black}{\gls{eh} is a process by which polymeric cellulose and hemicellulose molecules in biomass are broken down into monomeric sugars.} The default \gls{eh} setup is used with the exception of two input design parameters, the solid mass fraction of xylan (\textcolor{black}{hemicellulose}) and the \textcolor{black}{enzyme concentration (g/L)}, which are optimized by \gls{ac} to maximize the glucan (\textcolor{black}{sugars}) concentration. 

\begin{figure}
\centering
\includegraphics[width=1.0\columnwidth]{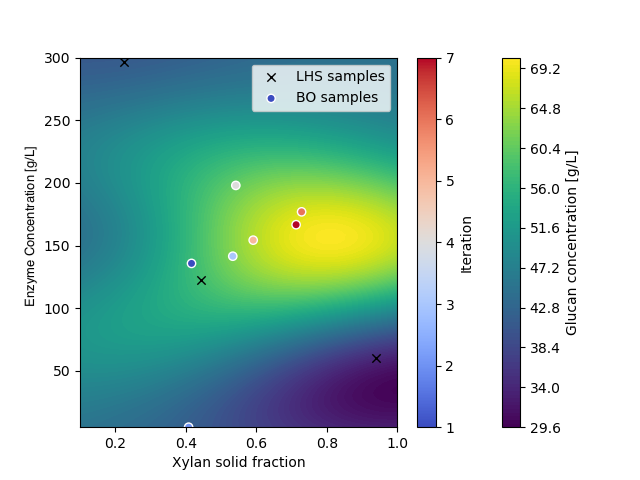}
\caption{Engineering design optimization for \gls{eh} reactor. The contour plot represents a Gaussian process model for glucan concentration trained on Vebio simulations indicated by the symbols. The black x's are three reactor designs with random samples using \gls{lhs}, and the colored circles are seven sequential samples using \gls{bo} with the \gls{ei} criteria for maximizing glucan concentration.}
\label{fig:ve}
\end{figure}

In Figure \ref{fig:ve}, the objective values are plotted versus these two simulation inputs. For the first three simulations (samples), the inputs are chosen randomly using \gls{lhs}. Then, \gls{ac} selects the next seven input pairs using \gls{bo}. The \gls{ei} acquisition function leverages the initial \gls{bo} samples to explore the parameter space and clusters the later samples where it expects the maximum to lie.

The benefit for the acquisition optimization is \gls{ei}. To choose the number of batches $I$, our heuristic is to make the batch wall clock budget $T_i$ equal to the wall-clock time of the longest simulation, to best approximate the performance of serial expected improvement, while permitting parallelization of low cost simulations. The number of remaining batches before the total wall-clock budget exhausts can be estimated by $T^r/T_i$. This implies the fraction of the resource budget that could be used if a second heuristic is introduced: use an equal fraction of the resource budget in each batch. Other strategies might use more of the budget on later batches as uncertainty is reduced.

\subsection{Autonomous Laboratory Processes Design: Guided Synthesis of GAN}

A key challenge in many scale-up problems is translating a laboratory proof of concept or prediction into a reliable, industrialized process. In material science, while advances in computing and \gls{ai} have enabled the prediction of new energy-relevant materials at an unprecedented rate, predicting synthesis pathways and reliably manufacturing such materials remain critical bottlenecks. 

The synthesis of a \gls{gan} thin film using \cgls{mme} encapsulates this difficulty. Though films have been grown with this method, controlling the vast arrays of associated process parameters (e.g. shutter timing, flux ratios, temperatures) and their influence on Ga droplet formation and film quality makes optimizing this process experimentally very expensive and difficult. Furthermore, while simulations can model this process and reproduce overall trends in, for example, defect discontinuity with shutter timing, there are substantial biases in the specific values in comparison with experiment \cite{gruber_2017}. This prohibits direct optimization of experimental parameters purely using simulation.

\subsubsection{AC formulation:}

\begin{figure*}
\centering
\includegraphics[width=1.0\textwidth]{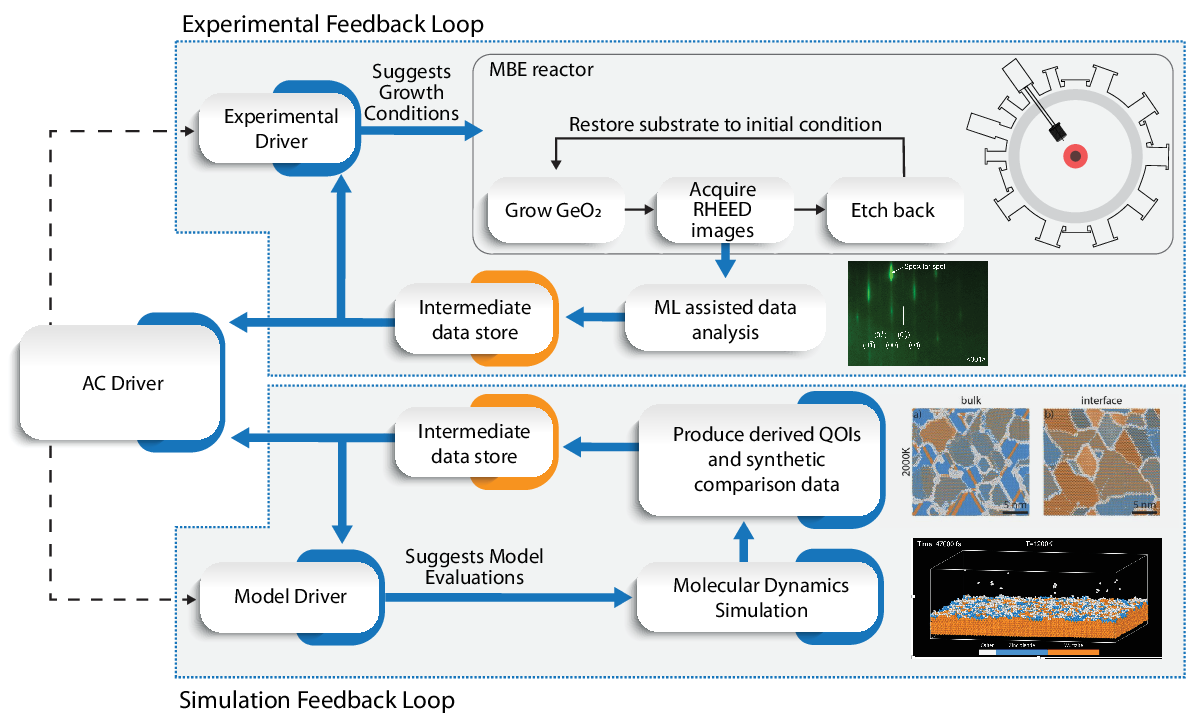}
\caption{The multi-fidelity experiment/simulation feedback loop for process control design. The experimental feedback loop iteratively grows and etches a thin film using suggested growth conditions from the experimental driver, storing in-situ data products. In parallel, the simulation feedback loop evaluates comparable molecular dynamic simulations and creates intermediate synthetic data products. The AC driver learns a model fidelity correction between the experiment and simulation and facilitates updating corresponding ML models.}
\label{fig:mat_synth}
\end{figure*}

We  cast this problem into the \gls{ac} framework by considering \gls{md} simulations and experimental results as separate fidelities \textcolor{black}{(low and high, respectively)}, with an outer loop problem of optimizing the process parameters for material quality. We then use the \gls{ac} framework to learn the underlying connection between the simulation and experiment, while taking into account experimental uncertainty with limited samples.

One essential consideration in using this approach with experimental infrastructure is maximizing experimental throughput and uptime while simultaneously utilizing the computational predictions and taking into account latency overhead and the disparate timescales of simulations and experiments. We address this within the \gls{ac} framework by extending the batch evaluation notion using a two feedback-loop drivers, as demonstrated in Figure \ref{fig:mat_synth}. This enables the experimental driver to query the AC driver as needed while being stable to communication interruptions and vast discrepancies in simulation/experiment timescale. The experimental feedback loop is enabled by an autonomous chemical re-etching procedure that allows a substrate to be returned to initial growth conditions repeatedly.

\subsection{Multi-solver Coupling: Adaptive Mesh and Algorithm Refinement for Perovskite Crystal Growth}

One common approach to solving multi-scale problems is to couple multiple solvers together, uni-directionally or bi-directionally, to resolve a large range of scales from the smallest, i.e., atomistic level, on the order of nanometers (molecules), to the largest, i.e., continuum level, on the order of meters (devices and reactors). These multi-solver simulation frameworks are particularly important for solving scale-up challenges because they enable accurate, predictive results that encompass the range of physics encountered in full-scale systems.

For illustration, at \gls{nrel}, we are using the \gls{amar} approach. It allows for active coupling between mesh-based \gls{pde} solvers at the continuum level along with particle-based solvers at the molecular level. The computational domain capturing the physics at the reactor and device levels is discretized using a background mesh. At locations of governing physical phenomena, e.g., catalyst/electrode surfaces, the mesh is refined using \gls{amr} and the continuum solver is replaced with a different physical approach (e.g., \gls{kmc} simulations), referred to as algorithm refinement. Information is transferred back and forth between the continuum and particle solvers in real-time. This has been used to model the dynamics inside a chemical vapor deposition reactor for perovskite crystal growth. 
The continuum flow solver transports 
species to a substrate. Several parallel \gls{kmc} simulations are then initiated for each location on the substrate with its corresponding gas phase compositions, pressure and temperature. The net flux of species at the surface is coupled to the continuum as a boundary condition.

\subsubsection{AC formulation:}
\begin{figure*}
\centering
\includegraphics[width=1.0\textwidth]{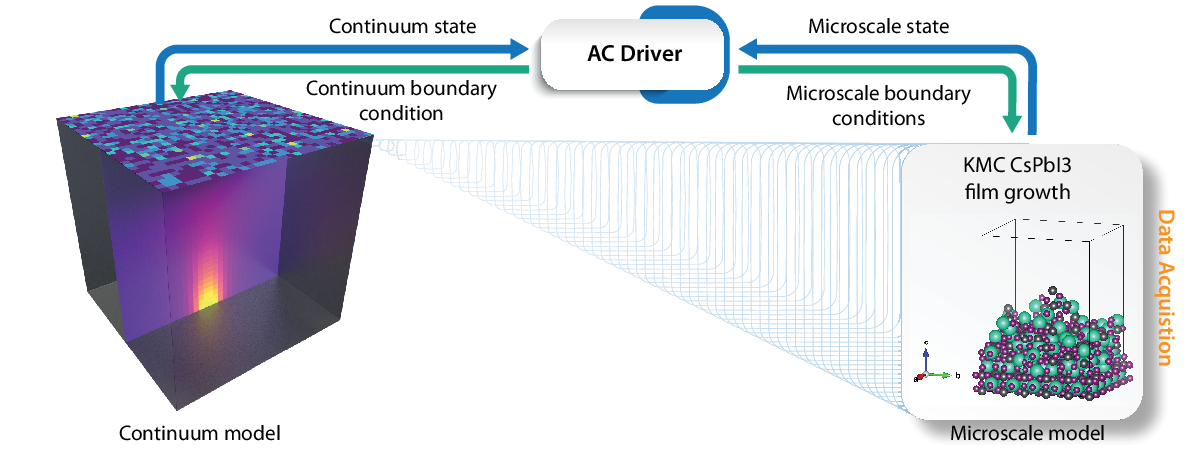}
\caption{A multi-species compressible flow impinges on a substrate (left). The perovskite crystal growth on the substrate is captured with a \gls{kmc} solver for adsorption, desorption, chemical reactions, and surface and bulk diffusion. The computed flux of species from the gas phase is passed back to the compressible flow code as a boundary condition for the next timestep. The \gls{ac} driver orchestrates the flow of information and the evaluation of the microscale model (either through a \gls{rom} or a high-fidelity \gls{kmc} evaluation).}
\label{fig:kmc_amar}
\end{figure*}

The outer loop problem in this case, as shown in Figure \ref{fig:kmc_amar}, 
is a macroscale continuum solver that launches microscale \gls{kmc} simulations for in-situ \gls{rom} generation. Instead of running a \gls{kmc} simulation at every grid cell at every time step to provide the surface boundary conditions for the continuum solver, \gls{ac} trains a surrogate model with a limited number of \gls{kmc} simulations and dynamically runs more simulations to reduce epistemic uncertainty. For the sampling strategy, $I$ is the number of time advancement steps in the macroscale solver. The resource and wall-clock budgets for the current step $B_i$ and $T_i$ are those allowed for the fraction of the each remaining quotas $B_r$ and $T_r$ to equal the fraction of time steps that remain to be completed. $M_\ell$ is the number of grid points in the macroscale solver because, at most, a microscale simulation could be run at every state observed in the macroscale simulation.
The benefit $r_{m\ell}$ is the variance of the surrogate model minus the variance of the next lowest level, if it exists.
This will prioritize simulations that are computationally inexpensive and correspond to large epistemic uncertainty at that operating condition on that fidelity level.

The \gls{ac} framework can also be used for engineering optimization for the design of the macro-scale conditions in the chemical vapor deposition reactor (e.g., ambient pressure, temperature, species composition, jet configuration), which determine the quality of the perovskite crystal grown at the substrate (e.g., number of defects, density, layer compositions, roughness).  
The \gls{ac} framework is also useful for managing the uncertainty arising from coupling a deterministic continuum model with a stochastic particle model. Tracking the effect of the resulting uncertainties on the \gls{qoi} (e.g., crystal layer growth, response of the flow at the substrate) is of vital importance.

\subsection{Site-specific Feedback Control: Flexible Building Electrical Loads}

We are expecting an unprecedented increase in the adoption of various clean energy technologies for buildings (e.g., solar panels, electric heat pumps, electric vehicle chargers) that are about to be deployed across millions of devices, geographies, and all types of end-use applications. 
This raises key questions around how to achieve grid resiliency objectives in the context of highly distributed control of buildings and their devices with very heterogeneous utilization schedules. 
Managing the expected building loads to facilitate their reduction and intelligent scheduling will enable communities to accomplish their decarbonization objectives while maintaining grid resiliency.
Furthermore, because of the large range of available building modeling tools across various fidelities and the large variability of building construction and utilization, buildings are a clear candidate to leverage the \Gls{ac} framework to automate the training of site-specific modeling and control.

\Gls{mpc} has proven to be a popular and successful method for the control of electrical loads and temperatures in buildings \cite{drgona_all_2020}. However, to obtain optimal control actions, \gls{mpc} requires a computationally expensive high-fidelity model, deeming \gls{mpc} computationally intractable for real-time implementation at scale.
One way to relax this computational burden is to use reduced-order models in the \gls{mpc} control loop. Though this solution may work well in the vicinity of training data, collecting enough relevant training data is expensive or even intractable for complex systems unless on-the-fly training is used. The \gls{ac} framework can leverage computing platforms across the edge, \gls{hpc}, and the cloud to train site-specific system models and device controllers to enable \gls{mpc} at scale.

\subsubsection{AC formulation:}
\begin{figure*}[hbt!]
\centering
\includegraphics[width=1.0\textwidth]{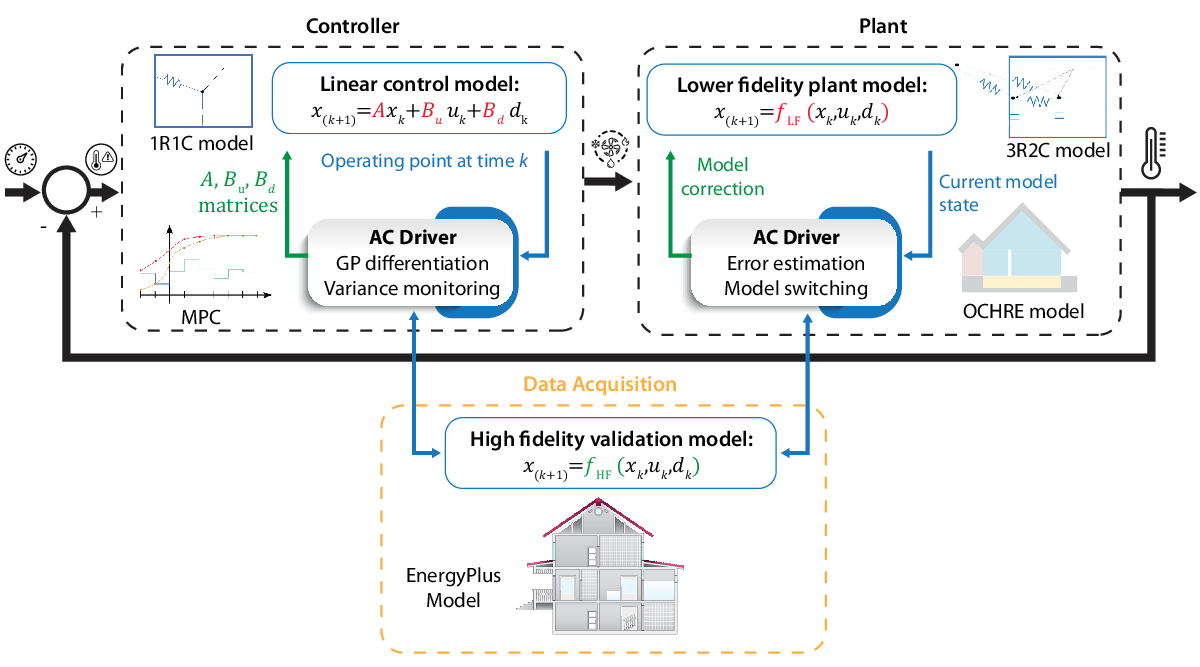}
\caption{The multi-fidelity building modelling and control methodology. The \gls{ac} driver (GP surrogate model) is used to adaptiveley update the low-fidelity linear control model based on the current conditions. In parallel, the plant model uses a second \gls{ac} driver to switch between models with various fidelities to balance accuracy and computational burden. The high fidelity EnergyPlus model updates these models with new data.}
\label{fig:grid}
\end{figure*}

As shown in Figure \ref{fig:grid}, \gls{ac} will play two distinct roles in dynamically updating the controller and the plant model of the system. 
The control model is based on a 1 resistor, 1 capacitor (1R1C) \gls{rom} with model coefficients predicted by \gls{ac}. A \gls{gp} surface maps the building states $x_k$, control inputs $u_k$, and disturbances $d_k$ at time $k$ to the state at time $k+1$ ($x_{k+1}$). The \gls{gp} surrogate model is then differentiated at the current conditions, returning the locally accurate \gls{rom}. This \gls{rom} is then used in the \gls{mpc}  loop to compute the optimal control actions $u_{k+1}$. An advantage of doing this in the \gls{ac} framework is that \gls{ac} also provides an estimate of the uncertainty at every operating point and can automatically launch a high-fidelity EnergyPlus simulation if the uncertainty exceeds a certain threshold, retrain the surrogate model, and achieve reduced uncertainty at the current operating point. This enables the \gls{gp} surrogate to be continuously updated with new data in an efficient manner.

\gls{ac} is also used to dynamically update the plant model from a fidelity hierarchy of three different models -- EnergyPlus, OCHRE, and 3R2C -- that describe the building thermodynamics with varying degrees of accuracy, and computational requirements. EnergyPlus, mentioned above, offers the highest accuracy at the highest computational cost. OCHRE introduces additional modeling assumptions to reduce computational cost. The 3R2C model is a low-cost gray-box model in which the building thermodynamics are abstracted to a resistance-capacitance circuit model, and the values of the resistors and capacitors are derived from historical time-series data. \gls{ac} orchestrates the use of these simulation paradigms to train a multi-fidelity model with maximum accuracy subject to wall-clock and computational resource budgets that are shared by all fidelities. The chosen acquisition function is the variance (uncertainty) of the multi-fidelity model. 

\section{Conclusions}

The unified \gls{ac} framework aims to derisk the scale-up of renewable energy technologies through uncertainty management for multifidelity simulations enabled by on-line learning from experiments and scientific computing subject to resource constraints. Consolidating these problems within a unified scheme for goal-setting and adaptively balancing computation cost, together with outlining concrete application examples, is an important step in accelerating impact and decreasing cross-disciplinary barriers.

Future work will include the full \gls{ac} framework core code release and subsequent core application demonstrations. Furthermore, we will continue to collaborate with domain experts to facilitate expansion into other application spaces in scale-up of renewable energy technologies. 
\section{Acknowledgements}
We would like to acknowledge the domain experts who supported the integration of \gls{ac} with our application areas.
Nick Wunder and Monte Lunacek developed the \gls{hero} framework and aided with the integration with \gls{ac}.  Davi Febba Marcelo, Stephen Schaefer, and Andriy Zakutayev developed the experimental infrastructure for the material synthesis portion, while Aidan Raver and Garritt Tucker provided simulation support.

This work was authored by the National Renewable Energy Laboratory (NREL), operated by Alliance for Sustainable Energy, LLC, for the U.S. Department of Energy (DOE) under Contract No. DE-AC36-08GO28308. This work was supported by the Laboratory Directed Research and Development (LDRD) Program at NREL. The views expressed in the article do not necessarily represent the views of the DOE or the U.S. Government. The U.S. Government retains and the publisher, by accepting the article for publication, acknowledges that the U.S. Government retains a nonexclusive, paid-up, irrevocable, worldwide license to publish or reproduce the published form of this work, or allow others to do so, for U.S. Government purposes. A portion of the research was performed using computational resources sponsored by the Department of Energy’s Office of Energy Efficiency and Renewable Energy and located at the National Renewable Energy Laboratory.  

\bibliographystyle{unsrt}
\bibliography{references} 

\begin{thebibliography}{10}

\bibitem{peherstorfer2018survey}
Benjamin Peherstorfer, Karen Willcox, and Max Gunzburger.
\newblock Survey of multifidelity methods in uncertainty propagation,
  inference, and optimization.
\newblock {\em Siam Review}, 60(3):550--591, 2018.

\bibitem{Kennedy2001}
Marc~C. Kennedy and Anthony O'Hagan.
\newblock Bayesian calibration of computer models.
\newblock {\em Journal of the Royal Statistical Society. Series B: Statistical
  Methodology}, 63:425--464, 2001.

\bibitem{meng2019}
Xuhui Meng and George~Em Karniadakis.
\newblock A composite neural network that learns from multi-fidelity data:
  Application to function approximation and inverse pde problems.
\newblock {\em Journal of Computational Physics}, 401:109020, 2020.

\bibitem{ghanem2017handbook}
Roger Ghanem, David Higdon, Houman Owhadi, et~al.
\newblock {\em Handbook of uncertainty quantification}, volume~6.
\newblock Springer New York, 2017.

\bibitem{hadjiconstantinou2006limits}
Nicolas~G Hadjiconstantinou.
\newblock The limits of navier-stokes theory and kinetic extensions for
  describing small-scale gaseous hydrodynamics.
\newblock {\em Physics of Fluids}, 18(11), 2006.

\bibitem{Jones1998}
Donald~R Jones, Matthias Schonlau, and William~J Welch.
\newblock Efficient global optimization of expensive black-box functions.
\newblock {\em Journal of Global Optimization}, 13:455--492, 1998.

\bibitem{shoemaker-paper}
Juliane M{\"{u}}ller and Christine~A. Shoemaker.
\newblock {Influence of ensemble surrogate models and sampling strategy on the
  solution quality of algorithms for computationally expensive black-box global
  optimization problems}.
\newblock {\em Journal of Global Optimization}, 60(2):123--144, 2014.

\bibitem{lam2016bayesian}
Remi Lam, Karen Willcox, and David~H Wolpert.
\newblock Bayesian optimization with a finite budget: An approximate dynamic
  programming approach.
\newblock {\em Advances in Neural Information Processing Systems}, 29, 2016.

\bibitem{shahriari2015taking}
Bobak Shahriari, Kevin Swersky, Ziyu Wang, Ryan~P Adams, and Nando De~Freitas.
\newblock Taking the human out of the loop: A review of bayesian optimization.
\newblock {\em Proceedings of the IEEE}, 104(1):148--175, 2015.

\bibitem{Young2023}
Ethan Young, Jonathan Stickel, Hariswaran Sitaraman, Andrew Glaws, and Jason
  Lischeske.
\newblock Virtual engineering (ve) (vebio).
\newblock [Computer Software] \url{https://doi.org/10.11578/dc.20220418.2}, Feb
  2022.

\bibitem{gruber_2017}
J.~Gruber, X.~W. Zhou, R.~E. Jones, S.~R. Lee, and G.~J. Tucker.
\newblock {Molecular dynamics studies of defect formation during
  heteroepitaxial growth of InGaN alloys on (0001) GaN surfaces}.
\newblock {\em Journal of Applied Physics}, 121(19):195301, 05 2017.

\bibitem{drgona_all_2020}
Ján Drgoňa, Javier Arroyo, Iago~Cupeiro Figueroa, David Blum, Krzysztof
  Arendt, Donghun Kim, Enric~Perarnau Ollé, Juraj Oravec, Michael Wetter,
  Draguna~L. Vrabie, and Lieve Helsen.
\newblock All you need to know about model predictive control for buildings.
\newblock {\em Annual Reviews in Control}, 50:190--232, 2020.

\end{thebibliography}

\section{Biographies} 

Rohit Chintala is a Senior Researcher at NREL at Golden, CO, 80401, USA. His research interests include building energy modeling, advanced control systems, and machine learning. Chintala received his Ph.D. in Mechanical Engineering from Texas A\&M University. Contact him at \href{rohit.chintala@nrel.gov}{rohit.chintala@nrel.gov}.

Marc Day is a Group Manager and Principal Scientist at \gls{nrel} in Golden, CO, 80401, USA. His research interests include turbulent reacting flows, high-performance computing and multi-scale/multi-fidelity modeling. Dr.\,Day received his Ph.D.\, in Applied Plasma Physics and Fusion Engineering from the University of California, Los Angeles. Contact him at \href{Marc.Day@nrel.gov}{Marc.Day@nrel.gov}.

Olga Doronina is a Data Scientist at \gls{nrel} at Golden, CO, 80401, USA. Her research interests include scientific machine learning, surrogate modeling and sensitivity analysis. Doronina received her Ph.D. in mechanical engineering from University of Colorado, Boulder. Contact her at \href{olga.doronina@nrel.gov}{olga.doronina@nrel.gov}.

Hilary Egan is a Computational Scientist at \gls{nrel} in Golden, CO, 80401, USA. Her research interests include scientific machine learning and integration of computation with experiments. Dr.\,Hilary Egan received her Ph.D.\, in Astrophysical and Planetary Science from the University of Colorado, Boulder. Contact at \href{hilary.egan@nrel.gov}{hilary.egan@nrel.gov}.

Kevin Patrick Griffin is Postdoctoral Research at \gls{nrel} in Golden, CO, 80401, USA. His research interests include data-driven modeling, turbulence, and \gls{cfd}. Dr.\,Griffin received his  Ph.D.\, in Mechanical Engineering from Stanford University. Contact him at \href{kevin.griffin@nrel.gov}{kevin.griffin@nrel.gov}.

Marc Henry de Frahan is a Computational Scientist at \gls{nrel} in Golden, CO, 80401, USA. His research interests include computational solvers for \glspl{pde}, high performance computing for heterogenous computing architectures, and multiphysics applications. Dr.\,Henry de Frahan received his Ph.D.\, in Mechanical Engineering from the University of Michigan, Ann Arbor. Contact him at \href{marc.henrydefrahan@nrel.gov}{marc.henrydefrahan@nrel.gov}.

Ryan King is a Senior Scientist at \gls{nrel} in Golden, CO, 80401, USA. His research interests include scientific machine learning, uncertainty quantification and turbulent flows. King received his Ph.D. in Mechanical Engineering from the University of Colorado, Boulder. Contact him at \href{ryan.king@nrel.gov}{ryan.king@nrel.gov}.

Ross Larsen is a Senior Scientist at \gls{nrel} in Golden, CO, 80401, USA. His research interests include electronic structure and dynamics in condensed phases, statistical mechanics of disordered systems, and molecular dynamics and Monte Carlo simulation methods. Larsen received his Ph.D. in Physics from Brown University. Contact him at \href{ross.larsen@nrel.gov}{ross.larsen@nrel.gov}.

Juliane Mueller is a Group Manager at NREL in Golden, CO, 80401, USA. Her research expertise includes derivative free optimization, surrogate modeling, and algorithm development for DOE relevant problems. Mueller received her Ph.D.  in applied math from Tampere University of Technology in Finland. She is a member of SIAM and INFORMS. Contact her at \href{Juliane.Mueller@nrel.gov}{Juliane.Mueller@nrel.gov}.

Jibonananda Sanyal is a Group Manager for Hybrid Energy Systems at NREL. His research interests include renewable energy, end-uses of energy, and applications of computing. Sanyal received his PhD in Computer Science from Mississippi State University. He is a Senior Member of IEEE. Contact him at \href{jibo.sanyal@nrel.gov}{jibo.sanyal@nrel.gov}.

Hariswaran Sitaraman is a computational scientist at NREL at Golden, CO, 80401, USA. His research interests include development of performance portable solvers and simulations of multiphase reacting and electrical driven flows. Sitaraman received his Ph.D. in Aerospace Engineering, University of Texas at Austin. Contact him at \href{hariswaran.sitaraman@nrel.gov}{hariswaran.sitaraman@nrel.gov}.

Deepthi Vaidhynathan is a Senior Researcher at NREL at Golden, CO, 80401, USA. Her research interests include energy system integration, grid modeling and simulation, and performance optimization for scientific software. Vaidhynathan received her M.S in Electrical Engineering from the University of Colorado, Boulder. She is a senior member at IEEE. Contact her at \href{Deepthi.Vaidhynathan@nrel.gov}{Deepthi.Vaidhynathan@nrel.gov}.

Dylan Wald is a PhD student at Colorado School of Mines in Golden, CO, 80401, USA. His research interests include optimal control, machine learning, and HVAC control in buildings. Contact him at \href{dylanwald@mines.edu}{dylanwald@mines.edu}.

Ethan Young is a Computational Scientist at NREL at Golden, CO, 80401, USA. His research interests include developing high-performance CFD tools for fluid-structure interaction problems, porous multiphase flows, and design and control optimization problems. Young received his Ph.D. in Biomedical Engineering from the University of California, Los Angles. Contact him at \href{ethan.young@nrel.gov}{ethan.young@nrel.gov}.


\end{document}